\documentclass[12pt]{amsart}
\usepackage{amscd,times,fullpage}
\usepackage[utf8x]{inputenc}

\usepackage{amsthm,amsmath}
\usepackage{amssymb}
\usepackage{graphicx}
\usepackage{amsthm}
\usepackage{enumerate}
\usepackage{dsfont}
\usepackage{tikz-cd}
\usetikzlibrary{cd}
\usetikzlibrary{positioning}
\usetikzlibrary{matrix}
\usepackage[hidelinks]{hyperref}

\hypersetup {
	pdfpagemode = {UseNone},
	pdftitle = {Approximation of Regular Sasakian Manifolds},
	pdfauthor = {Giovanni Placini},
	pdflang = {en-UK}
} 

\newcommand{\R}{\mathbb{R}}

\newcommand{\DD}{\mathbb{D}}       
\newcommand{\Z}{\mathbb{Z}}

\newcommand{\CP}{\mathbb{C}\mathrm{P}}
\newcommand{\N}{\mathbb{N}}

\newcommand{\Ric}{\mathrm{Ric}}

\newcommand{\C}{\mathbb{C}}            
\newcommand{\de}{\partial}

\newcommand{\K}{K\"{a}hler }

\newcommand{\OO}{\mathcal{O}}

\newcommand{\F}{\mathcal F}

\newcommand{\ov}[1]{\overline{#1}}

\newcommand{\deb}{\ov\partial}

\newcommand{\di}{{\operatorname{d}}}

\newcommand{\Id}{\operatorname{Id}}

\newcommand{\lra}{\longrightarrow}

\newcommand{\D}{\mathcal{D}}

\newcommand{\TM}{\text{T} M}

\newtheorem{theor}{Theorem}
\newtheorem{prop}[theor]{Proposition}
\newtheorem{defin}[theor]{Definition}

\newtheorem{ex}[theor]{Example}
\newtheorem{remark}[theor]{Remark}

\makeatletter
\newtheorem*{rep@theorem}{\rep@title}
\newcommand{\newreptheor}[2]{%
	\newenvironment{rep#1}[1]{%
		\def\rep@title{#2 \ref{##1}}%
		\begin{rep@theorem}}%
		{\end{rep@theorem}}}
\makeatother

\newreptheor{theor}{Theorem}

\usepackage{nicefrac}

\begin{document}

\title[Approximation of regular Sasakian manifolds]{Approximation of regular Sasakian manifolds}

	\author{G.~Placini}
\address{Dipartimento di Matematica e Informatica, Universit\'a degli studi di Cagliari, Via Ospedale 72, 09124 Cagliari, Italy}
\email{giovanni.placini@unica.it}

\date{\today ; {\copyright  G.~Placini 2023}}

	\thanks{The author is supported by INdAM and  GNSAGA - Gruppo Nazionale per le Strutture Algebriche, Geometriche e le loro Applicazioni, by GOACT - Funded by Fondazione di Sardegna and funded the National Recovery and Resilience Plan (NRRP), Mission 4 Component 2 Investment 1.5 - Call for tender No.3277 published on December 30, 2021 by the Italian Ministry of University and Research (MUR) funded by the European Union – NextGenerationEU. Project Code ECS0000038 – Project Title eINS Ecosystem of Innovation for Next Generation Sardinia – CUP F53C22000430001- Grant Assignment Decree No. 1056 adopted on June 23, 2022 by (MUR)}

\begin{abstract}
We investigate the problem of approximating a regular Sasakian structure by CR immersions in a standard sphere. 
Namely, we show that this is always possible for compact Sasakian manifolds.
Moreover, we prove an approximation result for non-compact $\eta$-Einstein manifolds via immersions in the infinite dimensional sphere. We complement this with several examples.
\end{abstract}

\maketitle

\section{Introduction and statements of the main results}\label{sectionint}

Sasakian geometry is often considered the odd dimensional analogue of \K geometry. This is due to the fact that a Sasakian manifold sits in a so-called "\K sandwich". Namely, a $2n+1$ dimensional Sasakian manifold comes with a \K $2n+2$ dimensional cone and a transverse \K geometry of dimension $2n$. 
This interplay translates to the fact that the solution of some problems in Sasakian geometry is equivalent to that of others in its older even dimensional analogue.
The problem considered in this paper falls into this case.
Namely, we ask whether a given regular Sasakian structure can be approximated by CR immersions in a standard sphere.
In analogy with a celebrated result of Tian, Ruan and Zelditch \cite{Tian90Approximation,Ruan98Approximation,Zelditch98Approximation}, it was proven in \cite{LoiPlacini22SasakiApproximations} that any compact Sasakian manifold is approximated by CR embeddings in a weighted sphere. 
Here we investigate two related questions.
Firstly, when the Sasakian structure is regular, it is natural to ask whether one can get a similar result to \cite[Theorem~1]{LoiPlacini22SasakiApproximations} under the requirement that the model space is a \textit{standard} Sasakian sphere.
Our first result shows that one can trade the injectivity of the embeddings for regularity in order to obtain immersions into the standard Sasakian sphere.

\begin{theor}\label{TheoMain1}
Let $(M,\eta,g)$ be a compact regular Sasakian manifold. Then there exist a sequence of CR immersions $\varphi_k:M\lra S^{2N+1}$ into standard Sasakian spheres such that suitable transverse homotheties of the induced structures converge to $(\eta, g)$ in the $C^\infty$-norm. 
\end{theor}
The solution to this problem is related to \K geometry in the following way. In the regular case the \K cone of a Sasakian manifold $M$ is the total space of a line bundle $L$ over a \K manifold $X$ (without the zero section).
It turns out that one can construct such immersions by means of an orthonormal basis for the space of sections of $L$ with respect to a certain scalar product. Finding such a basis is a classical problem in \K geometry deeply connected with the computation of Bergman kernels, special metrics and approximation of metrics, see for instance \cite{Donaldson01Embeddings1,Donaldson05Embeddings2,LoiMossa11BalancedBundles,Ma07Book,Tian90Approximation}.
Notice that one cannot avoid transverse homotheties because the transverse \K metric induced by an immersion in the standard sphere is an integral basic class while this is not necessarily true for the Sasakian structures we want to approximate.

Theorem~\ref{TheoMain1} above heavily relies on the compactness of the Sasakian manifold $M$.
Our second question is which conditions are sufficient for the existence of such approximation results in the noncompact case.
This is clearly a very broad question so we focus on the case of $\eta$-Einstein manifolds. 
In the compact case Cappelletti-Montano and Loi \cite{Cappelletti19EinsteinSpheres} studied immersions of compact regular $\eta$-Einstein manifold  into spheres with codimension $2$.
Here we prove an approximation result for (possibly noncompact) regular complete $\eta$-Einstein manifolds. 
\begin{theor}\label{TheoMain2}
Let $M$ be a complete regular $\eta$-Einstein manifold. Then the Sasakian structure on $M$ can be approximated by suitable $\D$-homotheties of a sequence of Sasakian structures induced by $CR$ embeddings in $S^\infty$.
\end{theor}
Also in this case, the immersions are constructed from a basis of the space of sections of a certain line bundle $L$ over a \K manifold $X$. 
The main difference with Theorem~\ref{TheoMain1} is the fact that $X$ is not compact so that the space of sections of $L$ could be infinite dimensional. 

As a particular case, all homogeneous Sasakian manifolds can be endowed with homogeneous $\eta$-Einstein metrics.
One should compare our result with \cite[Theorem~1.5]{LoiPlaciniZedda23SasakiHomogeneous} where the authors classify homogeneous Sasakian manifolds which admit an immersion in $S^\infty$.
In fact, the approximation in Theorem~\ref{TheoMain2} is constant for those homogeneous Sasakian manifolds which can be immersed in $S^\infty$.
In the last section of the we exhibit two examples of genuine approximations. 
Namely, we consider a Sasakian structure on $\C^*\times S^1$ and  $\DD^*\times S^1$ where $\C^*$ and disc $\DD^*$ are the punctured plane and disc respectively. 
We provide a sequence of embeddings of $\C^*\times S^1$ and  $\DD^*\times S^1$ into $S^\infty$ which approximate the given structures.
In terms of \K geometry, we compute the orthonormal basis of the space of sections of the trivial bundle over $\C^*$ and $\DD^*$.


\subsection*{Structure of the paper} The paper is organized as follows.  In  Section~\ref{sectionbackground} we review the basics of Sasakian geometry with particular focus on Sasakian immersions and regular Sasakian structures.
The remainder of the paper is divided into three sections. Namely, in Section~\ref{SecFinite}  and Section~\ref{SecInfinite} we prove Theorem~\ref{TheoMain1} and Theorem~\ref{TheoMain2} respectively.
Finally, Section~\ref{SecExamples} contains the computation of some explicit CR immersions of noncompact $\eta$-Einstein manifolds into $S^\infty$ approximating the given structure.

\section{Sasakian manifolds}\label{sectionbackground}

Sasakian geometry can be understood in terms of contact metric geometry and via the associated \K cone, cf. the monograph of Boyer and Galicki\cite{Boyer08Book}. We will present both formulations for the reader convenience, but we will focus mostly on the regular case for it is central in this paper.
In the following all manifolds and orbifolds are assumed to be connected.

A \textit{K-contact structure} $(\eta,\Phi,R,g)$ on a manifold $M$ consists of a contact form $\eta$ and an endomorphism 
$\Phi$ of the tangent bundle $\TM$ satisfying the following properties:
\begin{enumerate}
	\item[$\bullet$] $\Phi^2=-\Id+R\otimes\eta$ where $R$ is the Reeb vector field of $\eta$,
	\item[$\bullet$] $\Phi_{\vert\D}$ is an almost complex structure compatible with the symplectic form $\di\eta$ on $\D=\ker\eta$,
	\item[$\bullet$] the Reeb vector field $R$ is Killing with respect to the metric $g(\cdot,\cdot)=\dfrac{1}{2}\di\eta(\cdot,\Phi\cdot)+\eta(\cdot)\eta(\cdot)$.
\end{enumerate}   
Given such a structure one can consider the almost complex structure $I$ on the Riemannian cone $\big( M\times\R^+,t^2g+\di t^2\big)$ given by
\begin{enumerate}
	\item[$\bullet$] $I=\Phi$ on $\D=\ker\eta$, and
	\item[$\bullet$] $R=I(t\partial_t)_{\vert_{\{t=1\}}}$.
\end{enumerate} 
A \textit{Sasakian} \textit{structure} is a K-contact structure $(\eta,\Phi,R,g)$ such that the associated almost complex structure $J$ is integrable, and therefore $\left( M\times\R^+,t^2g+\di t^2,J\right)$ is K\"ahler. 
A \textit{Sasakian} manifold is a manifold $M$  equipped with a Sasakian structure $(\eta,\Phi,R,g)$.

Equivalently, one can define Sasakian manifolds in terms of \K cones. 
Namely, a Sasakian structure on a smooth manifold $M$ is defined to be a \K cone structure on $M\times\R^+=Y$.
That is, A \K structure $(g_Y,J)$ on $Y$ of the form 
$g_Y=t^2g+\di t^2$ where $t$ is the coordinate on $\R^+$ and $g$ a metric on $M$.
Then $(Y,g_Y,J)$ is called the \K cone of $M$ which, in turn, is identified with the submanifold $\{t=1\}$. The \K form on $Y$ is then given by
$$\Omega_Y=\dfrac{i}{2}\de\deb t^2.$$
The Reeb vector field on $Y$ is defined as 
$$R=J(t\de_t).$$
This defines a holomorphic Killing vector field with metric dual $1$-form
$$\eta=\dfrac{g_Y(R,\cdot)}{t^2}=\di^c\log t=i(\deb-\de)\log t$$
where $d^c=i(\deb-\de)$.
Notice that $J$ induces an endomorphism $\Phi$ of $\TM$ by setting 
\begin{enumerate}
	\item[$\bullet$] $\Phi=J$ on $\D=\ker\eta_{\vert_{\TM}}$, and
	\item[$\bullet$] $\Phi(R)=0$.
\end{enumerate} 
Equivalently, the endomorphism $\Phi$ is determined by $g$ and $\eta$ by setting
$$g(X,Z)=\dfrac{1}{2}\di\eta(X,\Phi Z)\ \ \mbox{ for }\ X,Z\in\D$$
It is easy to see that, when restricted to $M=\{t=1\}$, $(\eta,\Phi,R,g)$ is a Sasakian in the contact metric sense whose \K cone is $(Y,g_Y,J)$ itself. 
When this does not lead to confusion, we will use $R$ and $\eta$ to indicate both the objects on $Y$ and on $M$.

Since $g$ and $\eta$ are invariant for $R$, the Reeb foliation is transversally \K in the sense that the distribution $\D$ admits a \K structure $(g^T,\omega^T,J^T)$ which is invariant under $R$.
Explicitly, we have 
$$\omega^T=\dfrac{1}{2}\di\eta,\ \ J^T=\Phi_{\vert_{\D}}\ \mbox{ and } \ g^T(X,Z)=\dfrac{1}{2}\di\eta(X,J^T Z)=g_{\vert_{\D}}.$$
In particular, one can see that
\begin{equation}\label{EqTransvFormIsCurv}
	\omega^T=\dfrac{1}{2}\di\eta=\dfrac{i}{2}\di(\deb-\de)\log t=i\de\deb\log t.
\end{equation}

The Reeb vector field defines a foliation on $M$, called the Reeb foliation.
A very important dichotomy of Sasakian structures is given by the regularity of the leaves of the Reeb foliation.
Namely, if there exist foliated charts such that each leaf intersects a chart finitely many times, the Sasakian structure is called \textit{quasi-regular}.
Other wise it is called \textit{irregular}. Moreover, if every leaf intesects every chart only once, the sasakian structure is said to be \textit{regular}.
Compact regular and quasi-regular Sasakian manifold are fairly well understood due to the following result.

\begin{theor}[\cite{Boyer08Book}]\label{TheorStructure}
Let $(M,\eta,\Phi,R,g)$ be a quasi-regular  compact Sasakian manifold. Then the space of leaves of the Reeb foliation $(X,\omega,g_\omega)$ is a compact \K cyclic orbifold with integral \K form $\frac{1}{2\pi}\omega$ 
so that the projection $\pi:(M,g)\lra(X,g_\omega)$ is a Riemannian submersion. Moreover, $X$ is a smooth manifold if and only if the Sasakian structure on $M$ is regular. 

Viceversa, any principal $S^1$-orbibundle $M$ with Euler class $-\frac{1}{2\pi}[\omega]\in H^2_{orb}(X,\Z)$ over a compact \K cyclic orbifold $(X,\omega)$ admits a Sasakian structure.
\end{theor}
This result allows us to reformulate the geometry of a compact regular Sasakian manifold $M$ in terms of the  algebraic geometry of the \K manifold $X$. 
We will illustrate in detail this correspondence for its importance in the remainder of the paper. 
Let us first introduce the concept of $\D$-homothetic transformation of a Sasakian structure.

\begin{defin}[$\D$-homothety or a transverse homothety]\label{DefinTransverseHomothety}
Let $(M,\eta,\Phi,R,g)$ be a (not necessarily compact) Sasakian manifold and $a\in\R$ a positive number. One can define the Sasakian structure $(\eta_a,\Phi_a, R_a,g_a)$ from $(\eta,\Phi,R,g)$  as
$$ \eta_a=a\eta,\ \ \ \Phi_a=\Phi,\ \ \ R_a=\dfrac{R}{a},\ \ \  g_a=ag+(a^2-a)\eta\otimes\eta=ag^T+\eta_a\otimes\eta_a.$$

Equivalently, we can define the same structure on $M$ by setting a new coordinate on the \K cone as $\widetilde{t}=t^a$. It is clear from the formulation above that this induces on $M=\{\widetilde{t}=1\}=\{t=1\}$ the same Sasakian structure $(\eta_a,\Phi_a, R_a,g_a)$. We will call this structure the $\D_a$-homothety of $(\eta,\Phi,R,g)$.
\end{defin}

Now let the compact regular Sasakian manifold $(M,\eta,\Phi,R,g)$ be given and consider the projection $\pi:(M,g)\lra(X,\omega)$ given above. 
Notice that $\pi$ locally identifies the contact distribution $\D$ with the tangent space of $X$. Therefore,  up to $\D$-homothety, we have that $\pi^*(\omega)=\frac{1}{2}\di\eta$. 
Moreover, the endomorphism $\Phi$ determines the complex structure on $X$ and $g$ induces the \K metric $g_\omega$, i.e. $g^T=\pi^*g_\omega$.

In this case the class $\frac{1}{2\pi}[\omega]\in H^2(X,\Z)$ defines an ample line bundle $L$ over $X$.
Moreover, the cone $Y=M\times\R^+$ is identified with the complement of the zero section in $L^{-1}=L^*$ in the following way. Let $h$ be a hermitian metric on $L$ such that 
$$\omega=-i\de\deb\log h.$$
Then its dual $h^{-1}$ on $L^{-1}$ defines the second coordinate of $(p,t)\in M\times\R^+=L^{-1}\setminus\{0\}$ by
\begin{align}
	\nonumber t:L^{-1}\setminus\{0\} &\lra \R^+\\ 
	\label{EqMetricOnConeFromHermitian} (p,v)&\mapsto\vert v\vert_{h^{-1}_p}
\end{align}
where $v$ is a vector of $L^{-1}$ in the fiber over $p$.
With this notation the \K form on the \K cone $\left( M\times\R^+,t^2g+\di t^2,I\right)$ is given by

\begin{equation}
	\Omega=\dfrac{i}{2}\de\deb t^2.
\end{equation}
The Sasakian structure can be read from this data as
\begin{equation}\label{EqFromHermitianToSasakian}
	\omega^T=-i\de\deb\log h,\ \ \ \ \eta=i(\deb-\de)\log t.
\end{equation}

Therefore, the choice of a hermitian metric $h$ on an ample line bundle $L$ over a compact \K manifold $X$ completely determines a Sasakian structure on the $U(1)$-orbibundle associated to $L^{-1}$. 
The Sasakian manifold $(M,\eta,R,g,\Phi)$ so obtained is called a \textit{Boothby-Wang bundle} over $(X,\omega)$. 
Observe that, although the differentiable manifold is uniquely determined by $\frac{1}{2\pi}[\omega]$, the Sasakian structure does depend on, and is in fact determined by, the choice of $h$.

The most basic example is the standard Sasakian structure on $S^{2n+1}$, that is, the Boothby-Wang bundle determined by the Fubini-Study metric $h=h_{FS}$ on $\OO(1)$ over $\CP^n$. We give the details of this construction to further illustrate the formulation above.
\begin{ex}[Standard Sasakian sphere]\label{ExStandardSphere}\rm
Let $h=h_{FS}$ be the Fubini-Study hermitian metric on the holomorphic line bundle $\OO(1)$ over $\CP^n$. Recall that its dual metric $h^{-1}$ on $\OO(-1)\setminus\{0\}=\C^{n+1}\setminus\{0\}$ is given by the euclidean norm.
This defines a coordinate $t$ on the \K cone $\OO(-1)\setminus\{0\}=\C^{n+1}\setminus\{0\}=S^{2n+1}\times \R^+$. Namely, for coordinates $z=(z_0,z_1,\ldots,z_n)$ on $\C^{n+1}$ we have
\begin{align*}
	t:\C^{n+1} &\lra \R^+\\ 
	z&\mapsto\vert z\vert=\sqrt{\sum_{i=0}^{n}z_i\overline{z}_i}
\end{align*}
Now the \K metric on the cone is nothing but the flat metric
$$	\Omega_{flat}=\dfrac{i}{2}\de\deb t^2=\dfrac{i}{2}\sum\di z_i\wedge\di \overline{z}_i.$$
The Reeb vector field $R_0$ and the contact form $\eta_0$ read
$$	R_0=J(t\de_t)=i\sum z_i\de_{z_i}-\overline{z}_i\de_{\overline{z}_i},\ \ \ \eta_0=i(\deb-\de)\log t=\dfrac{i}{2t^2}\sum z_i\di\overline{z}_i-\overline{z}_i\di z_i.$$
It is clear that, when restricted to $S^{2n+1}$, $\eta_0$ and $R_0$, together with the round metric $g_0$ and the restriction $\Phi_0$ of $J$ to $\ker\eta_0$ give a Sasakian structure on $S^{2n+1}$.
This corresponds exactly to the Hopf bundle $S^{2n+1}\lra\CP^n$. Moreover, we have
$$\pi^*\omega_{FS}=\omega^T=\dfrac{1}{2}\di\eta_0=\dfrac{i}{2\vert z\vert^4}\sum_i\vert z_i\vert^2\di z_i\wedge\di \overline{z}_i-\sum_{i,j}\overline{z}_i z_j\di z_i\wedge\di \overline{z}_j$$
where $\pi:\C^{n+1}\setminus\{0\}\lra\CP^n$ is the standard projection.

Analogously, one can define the standard Sasakian structure on the infinite dimensional sphere $S^\infty=\{(z_o,z_1,\ldots)\in\ell^2(\C):\sum\vert z\vert^2=1\}$ (all sums are now infinite).
In this case the \K cone $S^\infty\times\R^+$ is the complex space $\ell^2(\C)\setminus\{0\}$ with the flat \K metric and the space of Reeb leaves is $\CP^\infty$.

\end{ex}
In general the space of leaves of the Reeb foliations $X$ is not even an orbifold. Nevertheless, when the Sasakian structure is regular and complete, $X$ is a \K manifold, see e.g. \cite{Reinhart59Foliated}.

We now switch our attention back to not necessarily compact Sasakian manifolds and recall another well known class of deformations of Sasakian structures, the so-called transverse \K transformations.
Namely, given a \K cone $Y=M\times\R^{+}$ we consider all \K metrics on $(Y,J)$ that are compatible with the Reeb field $R$. In other terms, these are potentials $\widetilde{t}^2$ such that $t\de_t=\widetilde{t}\de_{\widetilde{t}}$.
This means that the corresponding \K and contact forms satisfy 
$$\widetilde{\Omega}=\Omega+i\de\deb e^{2f},\ \ \ \widetilde{\eta}=\eta+d^cf$$
for a function $f$ invariant under $\de_t$ and $R$. Such functions are colled \textit{basic functions}.
We still need to identify the manifolds $\{\widetilde{t}=1\}$ and $\{t=1\}$.
This is done by means of the diffeomorphism
\begin{align*}
	F:Y&\lra Y\\
	(p,t)&\mapsto\left(p,te^{-f(p)}\right)
\end{align*}
which maps $\{t=1\}$ to $\{t=e^{-f(p)}\}=\{\widetilde{t}=1\}$. 
It is elementary to check that $\eta,\ R$ and $d^cf$ are invariant under $F$ so that $\widetilde{\eta}=\eta+d^cf$ holds on $M$.
Furthermore, the transverse \K forms are related by $\widetilde{\omega}^T=\omega^T+i\de\deb f$.
Notice that when the Sasaki structure is quasi-regular basic functions correspond to function on the base orbifold $X$. 
Thus, if $t$ comes from a hermitian metric $h^{-1}$ on $L^{-1}$, such a transformation is given by replacing $h^{-1}$ with $e^fh^{-1}$ for a function $f:X\lra\C$ such that $\omega+i\de\deb f>0$. 
This is equivalent to picking a different \K form $\widetilde{\omega}$ in the same class as $\omega$.
We summarise the above discussion in the following definition.
\begin{defin}[Transverse \K deformations]\label{DefinTransverseDeformation}
	Let  $(M,\eta,R,g,\Phi)$ be a Sasakian manifold with \K cone $(Y,J)$ and \K potential $t^2$. 
	A transverse \K transformation is given by replacing $t$ with $\widetilde{t}=e^{f}t$ for a basic function $f$ and leaving $(Y,J,R)$ unchanged. 
	When the Sasaki structure is quasi-regular and given as in \eqref{EqFromHermitianToSasakian}, a transverse \K transformation is given by replacing $h^{-1}$ with $e^fh^{-1}$.
\end{defin}

In this paper we are mostly interested with immersions and embeddings of Sasakian manifold. We recall the relevant definitions.
Two Sasakian manifolds $(M_1,\eta_1,R_1,g_1,\phi_1)$ and  $(M_2,\eta_2,R_2,g_2,\phi_2)$ are \text{equivalent} if there exists a diffeomorphism $\varphi\colon M_1\lra M_2$ such that
\begin{equation*}
	\varphi^*\eta_2=\eta_1\hspace{5mm}\mbox{and}\hspace{5mm}\varphi^*g_2=g_1.
\end{equation*}
If this holds then $\varphi$ also preserves the endomorphism $\phi_1$ and the Reeb vector field.
As implicitly intended in the definitions above, a Sasakian equivalence from a Sasakian manifold $(M,\eta,R,g,\phi)$ to itself is often called a \textit{Sasakian transformation} of $(M,\eta,R,g,\phi)$.

One can relax the condition on Sasakian equivalences to define Sasakian embedding and immersions. Namely, one does not request the map between Sasakian manifolds to be a diffeomorphism while requiring that it preserves the Sasakian structures.
In particular, given two Sasakian manifolds $(M_1,\eta_1,R_1,g_1,\phi_1)$ and  $(M_2,\eta_2,R_2,g_2,\phi_2)$, a \textit{Sasakian immersion (respectively embedding)} of $M_1$ in $M_2$ is an immersion (resp. embedding) $\varphi\colon M_1\lra M_2$ such that
\begin{align*}
	\varphi ^*\eta_2&=\eta_1,\hspace{20mm}\varphi ^*g_2=g_1,\\
	\varphi_*R_1&=R_2\hspace{5mm}\mbox{and}\hspace{5mm}\varphi_*\circ \phi_1=\phi_2\circ\varphi_*.
\end{align*}  
We can rephrase this definition in terms of the \K cone of the Sasakian manifolds $M_1$ and $M_2$. Namely, The map $\varphi$ satisfying the conditions above clearly extends to a map
\begin{align*}
	\widetilde \varphi\colon M_1\times\R&\lra M_2\times\R\\
	(p,t)&\mapsto (\varphi(p),t).
\end{align*}  
It is clear that if $\varphi$ is a Sasakian immersion (resp. embedding), then $\widetilde \varphi$ is a \K immersion (resp. embedding).

If, conversely, $Y_1$ and $Y_2$ are the \K cones of $M_1$ and $M_2$ with coordinates $t_1$ and $t_2$, then a \K immersion (resp. embedding) $\widetilde \varphi\colon Y_1 \lra Y_2$ such that $\widetilde\varphi^*(t_2)=t_1$ restricts to a Sasakian immersion (resp. embedding) $\varphi:M_1\lra M_2$.
Since it is often more useful to our purposes, we give the following

\begin{defin}[Sasakian immersion and embedding]\label{DefEmbedding}
	Let $M_1$ and $M_2$ be two Sasakian manifolds with \K cones $Y_1$ and $Y_2$ and coordinates $t_1$ and $t_2$ respectively. 
	A \textit{Sasakian immersion (respectively embedding)} of $M_1$ in $M_2$ is a \K immersion (resp. embedding) $\varphi\colon Y_1\lra Y_2$ such that $\varphi^*(t_2)=t_1$.
\end{defin}
\begin{remark}\rm
	Given the equivalence between a Sasakian immersion $M_1\lra M_2$ and a \K immersion of the \K cones,  with an abuse of notation, we will often denote both maps with the same letter. 
\end{remark}

A special class among Sasakian structures is that of $\eta$-Einstein structures. 
These are the Sasakian analogues of K\"ahler-Einstein metrics.
Namely, using the canonical splitting $TM=\D\oplus T_\F$ where $\D=\ker\eta$ and $T_\F$ denotes the tangent bundle to the Reeb foliation $\F$, write the metric as 
\begin{equation}\label{EqMetricSplit}
	g=g^T+\eta\otimes\eta.
\end{equation}

With an abuse of notation we write $g^T$ for both the transverse metric and the metric on $X$ in the quasi-regular case.
It follows from \eqref{EqMetricSplit} that the Riemannian properties of $M$ can be expressed in terms of those of the transverse K\"ahler geometry and of the contact form $\eta$. For instance, the Ricci tensor of $g$ is given by
\begin{equation}\label{EqRicciSplit}
	\Ric_g=\Ric_{g^T}-2g.
\end{equation}
A Sasakian manifold $(M,\eta,\phi,R,g)$ is said to be \textit{$\eta$-Einstein} if the Ricci tensor satisfies
\begin{equation}\label{EqEtaEinstein}
	\Ric_g=\lambda g+\nu \eta\otimes\eta
\end{equation}
for some constants $\lambda,\nu \in\R$. 
It follows from \eqref{EqRicciSplit} and \eqref{EqEtaEinstein} that a Sasakian manifold is $\eta$-Einstein with constants $(\lambda,\nu)$ if, and only if, its transverse geometry is K\"ahler-Einstein with Einstein constant $\lambda+2$.

	\subsection{CR immersions of regular and complete Sasakian manifolds into spheres}\label{SecGeneralSasakianEmbedding}
	
	We recall now some facts about CR immersions of Sasakian manifolds into finite and infinite dimensional standard spheres. 
	We only set the notation and report some useful results for us, the interested reader can refer to \cite[Section~5]{LoiPlaciniZedda23SasakiHomogeneous}

Let $M$ be a compact regular Sasakian manifold. 
By the Structure Theorem~\ref{TheorStructure}, $M$ is a $U(1)$-principal bundle $\pi:M\lra X$ over a compact \K manifold $(X,\omega)$ with $2\pi^*\omega=\di\eta$. 
Furthermore, $M$ is the unitary bundle associated to the line bundle $L^{-1}$ where $c_1(L)=[\omega]$.
This last condition implies that $L$ is ample. In other terms, $(X,L)$ is a polarised \K manifold.
Therefore, for $k\in\N$ large enough, the bundle $L^{\otimes k}=L^k$ is very ample, and we can define the Kodaira embedding $\psi_k:X\lra \CP^{N_k}$ where $\dim(H^0(L))=N_k+1$.
Then there exists a CR embedding $\varphi_k: M\lra S^{2N_k+1}$ of $M$ into the standard sphere covering the Kodaira embedding $\psi_k$ or, equivalently, a holomorphic embedding of $\varphi_k:Y\lra \C^{N_k+1}\setminus \{0\}$ of the \K cone $Y=M\times\R^+$ into the \K cone $S^{2N_k+1}\times\R^+$. In fact we have

\begin{prop}[{\cite[Proposition~5.1]{LoiPlaciniZedda23SasakiHomogeneous}}]\label{PropCompactEmbedding}
	Let $M$ be the compact regular Sasakian manifold determined by the Hermitian bundle $(L,h)$ over a compact projective manifold $X$. 
	For every integer $k>>0$ there exists a holomorphic embedding $\varphi_k:M\times\R^+\lra S^{2N_k+1}\times\R^+$ such that
	$\varphi_k^*(\tau)= B_kt^k$ where $B_k$ is the Bergman kernel of $L^k$, $\tau$ and $t$ are the coordinates on the second factor of $S^{2N_k+1}\times\R^+$ and $M\times\R^+$ respectively. 
\end{prop}

The same construction can be performed when the Sasakian manifold $M$ is the unitary bundle associated to the positive Hermitian bundle $(L,h)$ on a noncompact \K manifold $(X,\omega)$ with $\omega=-i\de\deb\log h$. In this case we cannot immerse $M$ into a finite dimensional sphere because the space of sections $H^0(L)$ is replaced by the Hilbert space $\mathcal{H}_{k,h}$ of integrable sections, see \cite{LoiPlaciniZedda23SasakiHomogeneous} for details. Nevertheless one gets the following noncompact analogue.

\begin{prop}\label{PropNoNCompactEmbedding}
	Let $M$ be the regular Sasakian manifold determined by the Hermitian bundle $(L,h)$ over a noncompact \K manifold $X$ and assume the space $\mathcal{H}_{k,h}$ is nontrivial. 
	Then there exists a holomorphic immersion $\varphi_k:M\times\R^+\lra S^{\infty}\times\R^+$ such that $\varphi_k^*(\tau)= \varepsilon_k t^k$ where $\varepsilon_k$ is the $\varepsilon$-function of $\mathcal{H}_{k,h}$, $\tau$ and $t$ are the coordinates on the second factor of $S^{\infty}\times\R^+$ and $M\times\R^+$ respectively. 
\end{prop}

\begin{remark}\rm\label{RmkPullbackCoordinate}
	Although $\varphi_k^*(h^{-1}_{FS})$ is not a hermitian metric on the line bundle $L^{-1}$ (it does not scale correctly under the $\C^*$-action), it defines a change of coordinate $(p,t)\mapsto (p,B_kt^k)$ (or $(p,t)\mapsto (p,\varepsilon_kt^k)$ in the noncompact case) on $M\times \R^+$ corresponding to the composition of the $\D_k$-homothetic transformation  ($t\mapsto t^k$) with a transverse \K deformation ($t^k\mapsto B_kt^k$).
\end{remark}

\section{Approximation of compact regular structures via immersions into spheres}\label{SecFinite}

\begin{proof}[Proof of Theorem~\ref{TheoMain1}]
Assume $(M,\eta,R,g,\Phi)$ to be a compact regular Sasakian manifold. 
Moreover, suppose we have performed a $\D$-homothetic transformation so that $M$ is the unit bundle $\pi:M\lra X$ associated to a holomorphic line bundle $L^{-1}$ over a projective manifold $(X,\omega)$ with $\pi^*\omega=\frac{1}{2}\di\eta$.

We can then apply Proposition~\ref{PropCompactEmbedding} to get a sequence of holomorphic immersions $\varphi_k:M\times\R^+\lra S^{2N_k+1}\times\R^+$ such that
$\varphi_k^*(\tau)= B_kt^k$ where $B_k$ is the Bergman kernel of $L^k$, $\tau$ and $t$ are the coordinates on the second factor of $S^{2N_k+1}\times\R^+$ and $M\times\R^+$ respectively. 
Notice that $\tau$ is the coordinate induced by the flat metric on $\C^{N_k+1}\setminus\{0\}=S^{2N_k+1}\times\R^+$ or, equivalently, by the hermitian metric $h_{FS}$ on $\OO(-1)$ whose curvature is $-\omega_{FS}$.

Now the $\frac{1}{k}$-transverse homothety of the structure induced on $M$ by the immersion into $S^{2N_k+1}$ is a transverse \K deformation of the original Sasakian structure determined by the Bergman kernel $B_k$, compare Definition~\ref{DefinTransverseDeformation} and Remark~\ref{RmkPullbackCoordinate}.
By \cite[Corollary~2]{Zelditch98Approximation} the first coefficient of the asymptotic expansion of the Bergman kernel $B_k$ smoothly converges to $1$ when $k$ goes to infinity. 
Therefore, the $\D_{\frac{1}{k}}$-homotheties of the structures determined by pullback coordinates $\varphi_k^*(\tau)= B_kt^k$ converge smoothly to $(\eta,R,g,\Phi)$.

We can resume the maps involved in the proof, with the notation of Section~\ref{SecGeneralSasakianEmbedding}, in the following diagram
$$
\begin{tikzcd}[column sep=large, row sep=large]
	(M, \eta_k,g_k) \arrow [d,"p_k"] \arrow [dd, bend right=50, "\pi"'] \arrow [dr,"\varphi_k"]\\		
	(M_k, \overline{\eta}_k,\overline{g}_k) \arrow [r,"\widetilde{\psi}_k"] \arrow[d,"\pi_k"] &
	(S^{2N_k+1},\eta_0,g_0) \arrow[d, "\pi_{FS}"] \\
	(X,\omega_k) \arrow[r,"\psi_k"] & \left(\CP^{N_k},\omega_{FS}\right)
\end{tikzcd}
$$
where $(M_k, \overline{\eta}_k,\overline{g}_k)$ is the unit bundle associated to $L^{-k}$ endowed with the Sasakian structure pulled back via $\psi_k$ and $(M, \eta_k,g_k)$ is the Sasakian structure determined by the coordinate $\varphi_k^*(\tau)= B_kt^k$.
\end{proof}

\begin{remark}\rm
Notice that we used a $\D$-homothety as the first step of the proof to get an actual Boothby-Wang bundle $\pi:M\lra X$. In order to avoid this and obtain the convergence to the original Sasakian metric, one can just compose the $\D_{\frac{1}{k}}$-homothety in our proof with the inverse of the homothetic transformation considered in the beginning.
\end{remark}

\section{Approximation of $\eta$-Einstein regular structures}\label{SecInfinite}

\begin{proof}[Proof of Theorem~\ref{TheoMain2}]
We cannot deduce that $M$ is an $S^1$-bundle over a K\"ahler manifold because $M$ is not necessarily compact.
Nevertheless, the Reeb foliation still defines a fibration $\pi\colon M\lra X$ over a K\"ahler manifold $(X,\omega)$ because $M$ is regular and complete, see \cite{Reinhart59Foliated}. 
Now the fibre of this fibration is either $\R$ or $S^1$.

Let us deal first with the case where the fibre is $S^1$.
Regardless of whether or not $M$ is compact, since the Sasakian structure on $M$ is regular and the fiber is $S^1$, it is the unit bundle of a line bundle $L^{-1}$ over $X$ such that $c_1(L)=[\omega]$.
Choose a hermitian metric $h$ on $L$ whose Ricci curvature form is $\omega$. 
Notice that $(X,\omega)$ is K\"ahler-Einstein because $M$ is $\eta$-Einstein.

We now invoke a result of Ma and Marinescu on the Bergman kernel of noncompact manifolds.
Namely, we apply \cite[Theorem~6.1.1]{Ma07Book} to the line bundle $L$ over $X$.
The hypotheses of this theorem are satisfied as $(X,\omega)$ is a K\"ahler-Einstein manifold so that there exists a positive constant $C$ such that $i\Ric(\omega)>C\omega$.
In our case this implies that the space of sections $\mathcal{H}_{k,h}$ is nontrivial so that Proposition~\ref{PropNoNCompactEmbedding} provides a sequence of holomorphic immersion $\varphi_k:M\times\R^+\lra S^{\infty}\times\R^+$ such that
$\varphi_k^*(\tau)= \varepsilon_k t^k$ where $\varepsilon_k$ is the $\varepsilon$-function of $\mathcal{H}_{k,h}$, $\tau$ and $t$ are the coordinates on the second factor of $S^{\infty}\times\R^+$ and $M\times\R^+$ respectively. 

Again by \cite[Theorem~6.1.1]{Ma07Book} (see also \cite{Arezzo13Kernel}) the $\varepsilon$-function $\varepsilon_k$ admits an asymptotic expansion whose first coefficient is $1$.
Therefore, taking the $\frac{1}{k}$-homothety of the Sasakian structure on $M$ defined by the pullback coordinate $\varphi_k^*(\tau)= \varepsilon_k t^k$ we get a sequence of structures which converge to the $\eta$-Einstein given one for $k\rightarrow\infty$.

Now the argument when the fiber is $\R$ easily follows from the previous one.
Namely, in this case the fibration is trivial, i.e. $M\cong X\times\R$.  
Since $\Z$ acts on $X\times\R$ by Sasakian isometries via the flow of the Reeb vector field, the quotient is the $\eta$-Einstein manifold $N=X\times S^1$ and the $\Z$-covering map $\widetilde \pi:M\lra N$ is a Sasakian immersion.
Now $N$ is an $\eta$-Einstein manifold fibering over a K\"ahler-Einstein manifold $X$ with fiber $S^1$.
By the previous case, there exists a sequence of CR immersions $\varphi_k:N\lra S^\infty$ such that suitable $\D$-homotheties of the induced structures converge on $N$ to the original $\eta$-Einstein structure. 
Therefore, the pullback to $M$ of such structures under $\widetilde\pi$ converge to the $\eta$-Einstein structure we began with.
Notice that these structures are transverse homotheties of the ones induced via the CR immersions $\widetilde \pi\circ\varphi_k:M\lra S^\infty$. That is, we can perform the transverse homotheties on $N$ or on $M$ interchangeably. 
This concludes the proof.
\end{proof}

\section{Explicit examples of approximations of $\eta$-Einstein structures}\label{SecExamples}

In this section we exploit the equivalence between polarizations $(L,h)$ of a \K manifold $X$ and Sasakian structures on a Boothby-Wang bundle over $X$ to describe explicitly some embeddings of noncompact inhomogeneous $\eta$-Einstein manifolds into $S^\infty$. Namely, we compute an orthonormal basis for the Hilbert space $\mathcal{H}_{k,h}$ of sections of a line bundle $L$ over a noncompact inhomogeneous \K manifold. This provides instances of approximations of inhomogeneous $\eta$-Einstein metrics which cannot be isometrically CR immersed in a sphere.

\begin{ex}[Fibring on the punctured plane]\rm
Consider the punctured plane $\C^*=\C\setminus\{0\}$ endowed with the complete Calabi-Yau metric $g^*_0$ induced by the \K form
\begin{equation}
	\omega^*_0=\frac{i}{2}\frac{\di z\wedge\di \overline{z}}{\vert z\vert^2}
\end{equation}
where $z$ is the coordinate on $\C^*$. 
Since this \K form admits a global potential $F=\frac{1}{2}\log^2\vert z\vert^2$, it is exact. 
Therefore, we can endow $\C^*\times S^1$, i.e. the unit bundle of the trivial bundle $L=\C^*\times\C$, with an $\eta$-Einstein structure. 
Namely, denoting the standard volume form on $S^1$ by $\alpha$, the contact form $\eta$ on $\C^*\times S^1$ is given by $\eta=\alpha+i(\deb-\de)\log F$. Moreover, the Sasakian metric is $g=g_0^*+\eta\otimes\eta$ and the endomorphism $\phi$ is given by the lift of the complex structure of $\C^*$ to the contact distribution.
We want to give an explicit expression of the embeddings of the $\eta$-Einstein manifolds $\C^*\times S^1$ just described into $S^\infty$. 

The \K space $(\C^*,\omega^*_0)$ and its polarizations were studied by Loi and Zuddas in \cite{Loi01Parma}. We report here the essential points which are relevant to our discussion. 
For any positive integer $k$
\begin{equation}
	h^k\left(f(z),f(z)\right)=e^{\frac{-k}{2}\log^2\vert z\vert^2 }\vert f(z)\vert^2
\end{equation}
is a hermitian metric on $L^k$ whose curvature is $k\omega_0^*$. 
By the discussion in the previous section, it is enough to compute an orthonormal basis of $\mathcal{H}_{k,h}$ to get the components of the embedding $\varphi_k$ of $L^{-k}\setminus\{0\}$ into $\ell^2(\C)$, cf. also \cite[Section~5]{LoiPlaciniZedda23SasakiHomogeneous}.
Namely, we need sections $s_j$ such that
\begin{equation}
	\langle s_j,s_j\rangle_k=\int_{\C^*}h^k\left(s_j(z),s_j(z)\right)\omega^*_0=\int_{\C^*}e^{\frac{-k}{2}\log^2\vert z\vert^2 }\vert s_j(z)\vert^2 \frac{i}{2}\frac{\di z\wedge\di \overline{z}}{\vert z\vert^2}=1
\end{equation}
and such that 	$\langle s_j,s_l\rangle_k=0$ for $j\neq l$.
It is easy to check that the functions $z^j$ for $j\in\Z$ are orthogonal and they form a basis of $\mathcal{H}_{k,h}$ for all $k$ since holomorphic functions are determined by their Laurent series.
A simple computation shows that
\begin{equation}
		\langle z^j,z^j\rangle_k=\int_{\C^*}e^{\frac{-k}{2}\log^2\vert z\vert^2 }\vert z\vert^{2j} \frac{i}{2}\frac{\di z\wedge\di \overline{z}}{\vert z\vert^2}=\sqrt{\frac{2}{k}}\pi^{\frac{3}{2}} e^{\frac{j^2}{2k}}.
\end{equation}
Hence an orthonormal basis for the Hilbert space $\mathcal{H}_{k,h}$ consists of the sections
$$s_{k,j}=\left(\frac{\sqrt{k}e^{\frac{-j^2}{2k}}}{\sqrt{2}\pi^{\frac{3}{2}}}\right)^{\frac{1}{2}}z^j$$
for $j\in\Z$.
In other words, the sections $s_{k,j}$ are the components of a holomorphic immersion $\varphi_k$ of $\C^*\times \C^*$ into $\ell^2(\C)$ and the potential of the induced transverse metric is 
$$F_k:=\varphi_k^*(\vert \cdot\vert^2)=\sum_{j\in\Z}\frac{\sqrt{k}e^{\frac{-j^2}{2k}}}{\sqrt{2}\pi^{\frac{3}{2}}}\vert z\vert^{2j}$$
so that the induced hermitian metric on $L^{-k}$ is 
$$\varphi_k^*(h_{FS}^{-1})\left(f(z),f(z)\right)=e^{F_k}\vert f(z)\vert^2.$$
One can check that the $k$-th root of this hermitian metric converges to (a multiple of) the metric $h=e^F\vert\cdot\vert^2$ without invoking \cite[Theorem~6.1.1]{Ma07Book}, see \cite[Theorem~3.6]{Loi01Parma} for a direct proof.
\end{ex}

\begin{ex}[Fibring on the punctured disc]\rm

As in the previous example we will construct noncompact Sasakian manifolds fibring over a noncompact inhomogeneous \K manifold $X$ with a global \K potential. 
Here we take $X=\DD^*=\{z\in\C:0<\vert z\vert^2<1\}$ equipped with the hyperbolic K\"ahler-Einstein metric 
\begin{equation}
	\omega^*_{hyp}=\frac{i}{2}\frac{\di z\wedge\di \overline{z}}{\vert z\vert^2\log^2\left(\vert z\vert^2\right)}
\end{equation}
whose potential is $F=-\log\left(-\log\vert z\vert^2\right)$. By Theorem~\ref{TheoMain2}
In analogy with the previous example, we can endow $\DD^*\times S^1$ with an $\eta$-Einstein structure with contact structure $\eta=\alpha+i(\deb-\de)\log F$. By Theorem~\ref{TheoMain2} this Sasakian structure can be  approximated by (suitable transverse homotheties) of structures induced by immersions of  $\C^*\times S^1$ into $S^\infty$. We want to give here the explicit expression of these immersions.

The \K space $(\C^*,\omega^*_{hyp})$ was studied in \cite{Ma21PuncturedAnnalen,Ma22PuncturedMathZ} in relation to Bergman kernels of punctured surfaces. The polarization we are interested is the $k$-th powers of the trivial line bundle endowed with the hermitian metric
\begin{equation}
	h^k\left(f(z),f(z)\right)=e^{k\log\left(-\log\vert z\vert^2\right) }\vert f(z)\vert^2.
\end{equation}
We compute an orthonormal basis of $\mathcal{H}_{k,h}$ to get the components of the embedding $\varphi_k$ of $L^{-k}\setminus\{0\}$ into $\ell^2(\C)$, cf. also \cite[Section~5]{LoiPlaciniZedda23SasakiHomogeneous}.
Namely, we need sections $s_j$ such that
\begin{equation}
	\langle s_j,s_j\rangle_k=\int_{\DD^*}h^k\left(s_j(z),s_j(z)\right)\omega^*_{hyp}=\int_{\DD^*}e^{k\log\left(-\log\vert z\vert^2\right) }\vert s_j(z)\vert^2 \frac{i}{2}\frac{\di z\wedge\di \overline{z}}{\vert z\vert^2\log^2\left(\vert z\vert^2\right)}=1
\end{equation}
and such that 	$\langle s_j,s_l\rangle_k=0$ for $j\neq l$.
It is easy to check that if a holomorphic function on $\DD^*$ has finite norm, then its Laurent expansion involves only the terms $z^j$ for positive $j\in\Z$.
The functions $z^j$ for $j>0$ are orthogonal and they form a basis of $\mathcal{H}_{k,h}$ for all $k$.
We can then compute
\begin{align*}
	\langle z^j,z^j\rangle_k&=\int_{\DD^*}e^{k\log\left(-\log\vert z\vert^2\right) }\vert z\vert^{2j} \frac{i}{2}\frac{\di z\wedge\di \overline{z}}{\vert z\vert^2\log^2\left(\vert z\vert^2\right)}= \frac{i}{2}\int_{\DD^*}\left(-\log\vert z\vert^2\right)^{k-2}\vert z\vert^{2j-2}\di z\wedge\di \overline{z}\\
	&=2\pi\int_{0}^{1}\left(-\log \rho^2\right)^{k-2}\rho^{2j-1}\di\rho
\end{align*}
where the last equality is obtained passing to polar coordinates. Substituting $e^x=\rho^2$ first and $-jx=w$ one gets
\begin{align*}
	\langle z^j,z^j\rangle_k&=	2\pi\int_{0}^{1}\left(-\log \rho^2\right)^{k-2}\rho^{2j-1}\di\rho=\pi\int_{-\infty}^{0}\left(-x\right)^{k-2}e^{jx}\di x\\
	&=\frac{\pi}{j^{k-1}}\int_{0}^{\infty}w^{k-2}e^{-w}\di w=\frac{\pi (k-2)!}{j^{k-1}}.
\end{align*}
Hence an orthonormal basis for the Hilbert space $\mathcal{H}_{k,h}$ consists of the sections
$$s_{k,j}=\left(\frac{j^{k-1}}{\pi(k-2)!}\right)^{\frac{1}{2}}z^j$$
for $j>0$ and these are the components of the holomorphic immersion $\varphi_k$ of $\DD^*\times \C^*$ into $\ell^2(\C)$.
In particular the potential of the induced transverse metric is 
$$F_k:=\varphi_k^*(\vert \cdot\vert^2)=\sum_{j>0}\frac{j^{k-1}\vert z\vert^{2j}}{\pi(k-2)!}$$
so that the induced hermitian metric on $L^{-k}$ is 
$$\varphi_k^*(h_{FS}^{-1})\left(f(z),f(z)\right)=e^{F_k}\vert f(z)\vert^2.$$
The $k$-th root of this hermitian metric converges to (a multiple of) the metric $h=e^F\vert\cdot\vert^2$ by \cite[Theorem~6.1.1]{Ma07Book}.
\end{ex}

\begin{remark}\rm
Notice that we can lift the $\eta$-Einstein structure of $\C^*\times S^1$ (respectively $\DD^*\times S^1$) to $\C^*\times\R$ (resp. $\DD^*\times \R$).
As in the proof of Theorem~\ref{TheoMain2}, by composing with the covering map, we can lift the immersions into $S^\infty$ too.

Observe none of these Sasakian manifolds are homogeneous Sasakian so that we provided explicit Sasakian immersions $\varphi_k$ of regular inhomogeneous $\eta$-Einstein manifolds into $S^\infty$ (when considered with the induced structure). 
This should be compared with \cite[Theorem~1.5]{LoiPlaciniZedda23SasakiHomogeneous} where it is proven that a homogeneous Sasakian manifolds can be immersed into $S^\infty$ if and only if its fundamental group is cyclic.
Our examples show that, if the manifold is not assumed to be homogeneous, there is no such restriction on the fundamental group.
\end{remark}

\bibliographystyle{amsplain}
		
\bibliography{biblio}

\providecommand{\bysame}{\leavevmode\hbox to3em{\hrulefill}\thinspace}
\providecommand{\MR}{\relax\ifhmode\unskip\space\fi MR }
\providecommand{\MRhref}[2]{%
  \href{http://www.ams.org/mathscinet-getitem?mr=#1}{#2}
}
\providecommand{\href}[2]{#2}
\begin{thebibliography}{10}

\bibitem{Arezzo13Kernel}
Claudio Arezzo, Andrea Loi, and Fabio Zuddas, \emph{Szeg\"{o} kernel, regular
  quantizations and spherical {CR}-structures}, Math. Z. \textbf{275} (2013),
  no.~3-4, 1207--1216. \MR{3127055}

\bibitem{Ma21PuncturedAnnalen}
Hugues Auvray, Xiaonan Ma, and George Marinescu, \emph{Bergman kernels on
  punctured {R}iemann surfaces}, Math. Ann. \textbf{379} (2021), no.~3-4,
  951--1002. \MR{4238257}

\bibitem{Ma22PuncturedMathZ}
\bysame, \emph{Quotient of {B}ergman kernels on punctured {R}iemann surfaces},
  Math. Z. \textbf{301} (2022), no.~3, 2339--2367. \MR{4437325}

\bibitem{Boyer08Book}
Charles~P. Boyer and Krzysztof Galicki, \emph{Sasakian geometry}, Oxford
  Mathematical Monographs, Oxford University Press, Oxford, 2008. \MR{2382957}

\bibitem{Cappelletti19EinsteinSpheres}
Beniamino Cappelletti-Montano and Andrea Loi, \emph{Einstein and {$\eta
  $}-{E}instein {S}asakian submanifolds in spheres}, Ann. Mat. Pura Appl. (4)
  \textbf{198} (2019), no.~6, 2195--2205. \MR{4031847}

\bibitem{Donaldson01Embeddings1}
S.~K. Donaldson, \emph{Scalar curvature and projective embeddings. {I}}, J.
  Differential Geom. \textbf{59} (2001), no.~3, 479--522. \MR{1916953}

\bibitem{Donaldson05Embeddings2}
\bysame, \emph{Scalar curvature and projective embeddings. {II}}, Q. J. Math.
  \textbf{56} (2005), no.~3, 345--356. \MR{2161248}

\bibitem{LoiPlacini22SasakiApproximations}
A.~Loi and G.~Placini, \emph{Any sasakian structure is approximated by
  embeddings into spheres}, 3 Oct 2022, preprint arXiv:2210.00790 [math.DG].

\bibitem{LoiPlaciniZedda23SasakiHomogeneous}
A.~Loi, G.~Placini, and M.~Zedda, \emph{Immersions into {S}asakian space
  forms}, 9 May 2023, preprint arXiv:2305.05509 [math.DG].

\bibitem{Loi01Parma}
A.~Loi and D.~Zuddas, \emph{Some remarks on {B}ergmann metrics}, Riv. Mat.
  Univ. Parma (6) \textbf{4} (2001), 71--86. \MR{1878012}

\bibitem{LoiMossa11BalancedBundles}
Andrea Loi and Roberto Mossa, \emph{Uniqueness of balanced metrics on
  holomorphic vector bundles}, J. Geom. Phys. \textbf{61} (2011), no.~1,
  312--316. \MR{2747002}

\bibitem{Ma07Book}
Xiaonan Ma and George Marinescu, \emph{Holomorphic {M}orse inequalities and
  {B}ergman kernels}, Progress in Mathematics, vol. 254, Birkh\"{a}user Verlag,
  Basel, 2007. \MR{2339952}

\bibitem{Reinhart59Foliated}
Bruce~L. Reinhart, \emph{Foliated manifolds with bundle-like metrics}, Ann. of
  Math. (2) \textbf{69} (1959), 119--132. \MR{107279}

\bibitem{Ruan98Approximation}
Wei-Dong Ruan, \emph{Canonical coordinates and {B}ergmann [{B}ergman] metrics},
  Comm. Anal. Geom. \textbf{6} (1998), no.~3, 589--631. \MR{1638878}

\bibitem{Tian90Approximation}
Gang Tian, \emph{On a set of polarized {K}\"{a}hler metrics on algebraic
  manifolds}, J. Differential Geom. \textbf{32} (1990), no.~1, 99--130.
  \MR{1064867}

\bibitem{Zelditch98Approximation}
Steve Zelditch, \emph{Szeg\={o} kernels and a theorem of {T}ian}, Internat.
  Math. Res. Notices (1998), no.~6, 317--331. \MR{1616718}

\end{thebibliography}
		
\end{document}